\documentclass{article}
\usepackage{amsmath, amssymb, latexsym, amscd}
\usepackage{amsfonts}
\newtheorem{definition}{Definition}[section]
\newtheorem{example}[definition]{Example}
\newtheorem{lemma}[definition]{Lemma}
\newtheorem{theorem}[definition]{Theorem}
\newtheorem{corollary}[definition]{Corollary}
\newtheorem{proposition}[definition]{Proposition}
\newtheorem{conjecture}[definition]{Conjecture}
\newtheorem{remark}[definition]{Remark}
\newenvironment{proof}{\noindent {\it Proof}:}{}

\author{Christopher A. Francisco\\Department of Mathematics, University of Missouri\\Columbia, MO 65211-4100\\E-mail: chrisf@math.missouri.edu}
\title{New approaches to bounding the multiplicity of an ideal}
\date{}
\begin{document}
\maketitle

\begin{abstract}
We use the theory of resolutions for a given Hilbert function to investigate the multiplicity conjectures of Huneke and Srinivasan and Herzog and Srinivasan. To prove the conjectures for all modules with a particular Hilbert function, we show that it is enough to prove the statements only for elements at the bottom of the partially ordered set of resolutions with that Hilbert function. This enables us to test the conjectured upper bound for the multiplicity efficiently with the computer algebra system Macaulay 2 \cite{M2}, and we verify the upper bound for many Artinian modules in three variables with small socle degree. Moreover, with this approach, we show that though numerical techniques have been sufficient in several of the known special cases, they are insufficient to prove the conjectures in general. Finally, we apply a result of Herzog and Srinivasan on ideals with a quasipure resolution to prove the upper bound for Cohen-Macaulay quotients by ideals with generators in high degrees relative to the regularity.
\end{abstract}

\section{Introduction}\label{intro}

Conjectures of Huneke and Srinivasan and Herzog and Srinivasan in the 1990s have led to considerable effort recently to bound the multiplicity of a homogeneous ideal in terms of the shifts in its graded free resolution. In this paper, we introduce new approaches to these conjectures to prove some special cases and gain some insight into why it has been so difficult to make general progress. Throughout, $R=k[x_1,\dots,x_n]$, the polynomial ring in $n$ variables over a field $k$. All ideals are homogeneous.

The motivation for comparing the multiplicity of a module to products of the shifts in the graded free resolution comes from a paper of Huneke and Miller \cite{HunekeMiller}. They focused on Cohen-Macaulay graded modules $R/I$ with a pure resolution, meaning that all the minimal generators of $I$ have the same degree $d_1$, all the minimal first syzygies of $I$ have the same degree $d_2$, etc. Huneke and Miller proved:

\begin{theorem}
\label{pureres}
\emph{(Huneke-Miller)} Let $I \subset R$ be a homogeneous ideal of codimension $c$ such that $R/I$ is Cohen-Macaulay with minimal graded free resolution \[ 0 \rightarrow R^{\beta_c}(-d_c) \rightarrow \cdots \rightarrow R^{\beta_1}(-d_1) \rightarrow R \rightarrow R/I \rightarrow 0 \] and multiplicity $e(R/I)$. Then \[ e(R/I) = \frac{1}{c!} \prod_{i=1}^c d_i.\]
\end{theorem}

As Huneke and Miller point out, because $e(R/I)$ is an integer, Theorem~\ref{pureres} places strong conditions on the shifts that can occur for a module with a pure resolution. It is natural to ask what the analagous statement would be for graded modules $R/I$ that do not have a pure resolution. For each $i$, let \[ m_i= \min \{j : \beta_{ij}(R/I) \not = 0 \} \quad \hbox{and} \quad M_i= \max \{j : \beta_{ij}(R/I) \not = 0 \},\] where the $\beta_{ij}(R/I)$ are the graded Betti numbers of $R/I$. Thus $m_i$ is the minimal degree of a syzygy at step $i$ of the minimal graded free resolution of $R/I$, and $M_i$ is the corresponding maximum. The goal is to bound $e(R/I)$ in terms of the $m_i$ and $M_i$.

Huneke and Srinivasan made the next conjecture in the Cohen-Macaulay case, and Herzog and Srinivasan generalized it to the case in which $R/I$ is not Cohen-Macaulay.

\begin{conjecture}
\label{mainconj}
\emph{(Huneke-Srinivasan, Herzog-Srinivasan)} Let $I \subset R$ be a homogeneous ideal of codimension $c$. If $R/I$ is Cohen-Macaulay, then \[ \frac{1}{c!} \prod_{i=1}^c m_i \le e(R/I) \le \prod_{i=1}^c M_i. \] If $R/I$ is not Cohen-Macaulay, then only \[ e(R/I) \le \frac{1}{c!} \prod_{i=1}^c M_i. \]
\end{conjecture}

Thus the conjecture of Herzog and Srinivasan is that the Cohen-Macaulay hypothesis is irrelevant for the upper bound. In the non-Cohen-Macaulay case, the codimension is less than the projective dimension, so the upper bound would only depend on the first $c$ steps of the resolution and hence can be considered a stronger statement than in the Cohen-Macaulay case. Note that the lower bound fails badly if $R/I$ is not Cohen-Macaulay; for example, if $I=(a^2,ab) \subset S=k[a,b]$, then $e(S/I)=1$, codim $I=1$, and $m_1=2$, but $2 \not < 1$.

There is a growing body of literature proving special cases of Conjecture~\ref{mainconj}. We mention a number of cases here; for a detailed discussion of what is known about the conjectures, see the expository paper \cite{FranciscoSrinExpos}. In codimension two, the conjectures are completely solved. Herzog and Srinivasan did the Cohen-Macaulay case in \cite{HerzogSrin}, and R\"{o}mer proved the non-Cohen-Macaulay case in \cite{Roemer}. Gold also has results on codimension two lattice ideals \cite{Gold}. Recently, Migliore, Nagel, and R\"{o}mer proved stronger bounds for Cohen-Macaulay codimension two ideals, showing that the bounds are sharp if and only if the ideals have pure resolutions \cite{MNR:mult}. They also proved the Gorenstein codimension three case, generalizing earlier results of Herzog and Srinivasan, and gave the same sharpness result as in the codimension two case. Further, in the more recent paper \cite{MNR2}, Migliore, Nagel, and R\"{o}mer found a lower bound for the degree of non-arithmetically Cohen-Macaulay curves in $\mathbb P^3$ and proved some cases of a generalization of Conjecture~\ref{mainconj} for modules. They also did an extensive analysis of the behavior of Conjecture~\ref{mainconj} under basic double G-linkage, proving the conjecture for standard determinantal ideals, a result obtained independently by Mir\'o-Roig in \cite{M-R} and then Herzog and Zheng in \cite{HZ}. Herzog and Zheng also proved that, in most of the special cases of the conjecture that are known, the bounds of Conjecture~\ref{mainconj} are sharp if and only if $R/I$ is Cohen-Macaulay and has a pure resolution, and they explored powers of ideals that are known to satisfy the bounds \cite{HZ}. In \cite{HerzogSrin}, Herzog and Srinivasan proved Conjecture~\ref{mainconj} for Cohen-Macaulay modules with a quasipure resolution, complete intersections, and stable ideals (among other cases), and R\"{o}mer used the result for stable ideals to prove the conjectures for the more general class of componentwise linear ideals \cite{Roemer}. Additionally, Guardo and Van Tuyl proved the conjectures for powers of a complete intersection in \cite{GuardoVT}. On the geometric side, Gold, Schenck, and Srinivasan proved Conjecture~\ref{mainconj} for some configurations of points in projective space in \cite{GSS}, and the conjecture is known for small sets of general fat points in $\mathbb P^n$ as well \cite{FranciscoFat}. Finally, Herzog and Srinivasan have some results on weaker upper bounds than in Conjecture~\ref{mainconj} using the Taylor resolution \cite{HerzogSrinMons}. 

Despite this work, we have made little general progress on Conjecture~\ref{mainconj}. One of the aims of this paper is to give some insight into why the conjectures have been so difficult. We do this by relating Conjecture~\ref{mainconj} to a question in the theory of resolutions for a given Hilbert function. This connection also gives us a sufficient condition with which we can efficiently test the conjectures for modules with a fixed Hilbert function.

Our paper is organized as follows. In Section~\ref{reshf}, we review some results about the possible graded free resolutions that can occur for a module with a given Hilbert function. This theory allows us to reduce Conjecture~\ref{mainconj} to considering a finite number of sets of graded Betti numbers for each Hilbert function. We present some computational work based on this reduction in the next section, using Macaulay 2 to verify the upper bound of Conjecture~\ref{mainconj} for a large number of Artinian ideals in three variables. As a result of the computations we made, we show that although numerical techniques have sufficed to prove some special cases of Conjecture~\ref{mainconj}, they are not enough in general, and in Remark~\ref{implications}, we discuss implications for the interplay between the multiplicity conjectures and problems on resolutions for a given Hilbert function. In the final section, we discuss a technique for attacking the upper bound in the Cohen-Macaulay case, using truncation to eliminate superfluous information in the graded free resolution. We use trunctation and Herzog and Srinivasan's result in the Cohen-Macaulay quasipure case to prove the upper bound for ideals with generators in high degree.

This work is part of the author's Ph.D. thesis at Cornell University. I thank the Cornell Graduate School for its fellowship support, Mike Stillman for his help with the results in this paper, and Juan Migliore for his insightful comments. I also gratefully acknowledge Grayson and Stillman's computer algebra system Macaulay 2 \cite{M2} on which I did the computations in this paper.


\section{Resolutions for a given Hilbert function}\label{reshf}

In his seminal 1890 paper \cite{Hilbert}, Hilbert used the graded free resolution of a module to compute the Hilbert function. Given any graded free resolution of $R/I$, we know precisely what the Hilbert function of $R/I$ is. The converse question has much more substance: Given a Hilbert function, what are the possible minimal graded free resolutions (or sets of graded Betti numbers) that occur for modules with that Hilbert function?

To study this question, we impose a partial order on resolutions of modules with a given Hilbert function. Suppose $R/I$ and $R/J$ have the same Hilbert function, and we wish to compare their graded Betti numbers. We say that $\beta^{R/I} \le \beta^{R/J}$ if and only if $\beta^{R/I}_{ij} \le \beta^{R/J}_{ij}$ for all $i$ and $j$. This is a strong conditon; in particular, there are likely to be a number of incomparable resolutions for modules with a particular Hilbert function. 

It is natural to ask about the structure of the resulting partially ordered set of resolutions for a given Hilbert function. The behavior at the top of the partial order is particularly nice. Recall that an ideal $L$ is called a lexicographic ideal if it is a monomial ideal generated in each degree $d$ by the first $\dim_k L_d$ monomials in descending lexicographic order. (For example, $(a^2,ab,ac^2,b^4,b^3c)$ is a lexicographic ideal in $k[a,b,c]$.) Given a homogeneous ideal $I \subset R$, there is always a lexicographic ideal in $R$ with the same Hilbert function. The following result, due independently to Bigatti \cite{Bigatti} and Hulett \cite{Hulett} in characteristic zero and Pardue \cite{Pardue} in positive characteristic, shows that the resolution of the lexicographic ideal is always the unique top element in the partial order.

\begin{theorem}
\label{bhp}
\emph{(Bigatti, Hulett, Pardue)} Let $I \subset R$ be a homogenous ideal, and let $L$ be the lexicographic ideal with the same Hilbert function. Then $\beta^{R/I} \le \beta^{R/L}$.
\end{theorem}

Theorem~\ref{bhp} makes the search for all possible sets of graded Betti numbers for a given Hilbert function into a finite problem since all resolutions must lie under that of the lexicographic ideal. We can say a bit more: Let \{$\beta_{ij}$\} be a set of graded Betti numbers. Fix some integers $r$ and $j$, and replace $\beta_{rj}$ by $\beta_{rj}-1$ and $\beta_{r+1,j}$ by $\beta_{r+1,j}-1$. We call this a \textbf{consecutive cancellation}. A theorem of Peeva shows that all sets of graded Betti numbers for modules with a given Hilbert function are obtained by a sequence of consecutive cancellations in the resolution of the lexicographic ideal \cite{Peeva}. Thus when we speak of a potential Betti diagram in this paper, we mean a Betti diagram obtained from consecutive cancellations in the Betti diagram of a lexicographic ideal. Peeva's result places even stronger restrictions than Theorem~\ref{bhp} on what configurations of graded Betti numbers can occur. Note, however, that at the bottom of the partial order, the situation is not as simple as it is at the top. Charalambous and Evans have shown that there may be incomparable minimal elements in the partial order for a particular Hilbert function \cite{CE:HF}. The structure of the partially ordered set is an active area of research; see, for example, the work of Richert \cite{Richert}.

How does this relate to Conjecture~\ref{mainconj}? First, because lexicographic ideals are stable, by a result of Herzog and Srinivasan \cite{HerzogSrin}, they satisfy Conjecture~\ref{mainconj}. Thus the conjectured bounds hold for the resolution at the top of the partial order. What about the other elements in the partial order?

To consider this question, let $I$ be a homogeneous ideal of codimension $c$ in $R=k[x_1,\dots,x_n]$. Suppose $R/I$ satisfies the bounds of Conjecture~\ref{mainconj}. (If $R/I$ is not Cohen-Macaulay, then by this assumption we mean that $R/I$ satisfies the upper bound, and one makes the obvious adjustments in the discussion that follows.) Then, if $m_i^{R/I}$ and $M_i^{R/I}$ are the minimum and maximum shifts at each step of the minimal resolution of $R/I$ as in Section 1, we have \[ \frac{1}{c!} \prod_{i=1}^c m_i^{R/I} \le e(R/I) \le \frac{1}{c!} \prod_{i=1}^c M_i^{R/I}. \] 

Let $J$ be a homogeneous ideal in $R$ satisfying the following conditions: $R/J$ has the same Hilbert function as $R/I$, and $\beta^{R/I} \le \beta^{R/J}$. (Since $I$ and $J$ are both ideals in $R$ with the same Hilbert series, they have the same codimension.) Then because $\beta_{ij}^{R/I} \le \beta_{ij}^{R/J}$ for all $i$ and $j$, \[ m_i^{R/J} \le m_i^{R/I} \quad \hbox{and} \quad M_i^{R/I} \le M_i^{R/J}; \] the resolution of $R/J$ has all the terms of the resolution of $R/I$ plus possibly more, so the minimum shifts can be lower, and the maximum shifts can be higher. In addition, because $I$ and $J$ have the same Hilbert function, and the Hilbert function determines the multiplicity, $e(R/I)=e(R/J)$. Putting these facts together, we have \[ \frac{1}{c!} \prod_{i=1}^c m_i^{R/J} \le \frac{1}{c!} \prod_{i=1}^c m_i^{R/I} \le e(R/I)=e(R/J) \le \frac{1}{c!} \prod_{i=1}^c M_i^{R/I} \le \frac{1}{c!} \prod_{i=1}^c M_i^{R/J}. \]

Therefore $R/J$ also satisfies the bounds of Conjecture~\ref{mainconj}. Consequently, we have proven the following proposition.

\begin{proposition}
\label{multsamehf}
Let $I$ and $J$ be homogeneous ideals in $R=k[x_1,\dots,x_n]$. Suppose that $I$ and $J$ have the same Hilbert function, $\beta^{R/I} \le$ $\beta^{R/J}$, and $R/I$ satisfies the bounds in Conjecture~\ref{mainconj}. Then Conjecture~\ref{mainconj} also holds for $R/J$. \hfill $\square$
\end{proposition}

The condition that $\beta^{R/I} \le \beta^{R/J}$ in Proposition~\ref{multsamehf} is actually a bit stronger than what we need. It is enough, for example, to have that $\beta_{ij}^{R/I} \not =0$ implies that $\beta_{ij}^{R/J} \not =0$. 

Proposition~\ref{multsamehf} reduces checking Conjecture~\ref{mainconj} for modules with a given Hilbert function to verifying the bound for a finite number of sets of graded Betti numbers. To see this, fix a Hilbert function. Suppose we can determine all the resolutions that are minimal in the partial order for that Hilbert function. If we can show that Conjecture~\ref{mainconj} holds for all the minimal elements in the partial order, we can use Proposition~\ref{multsamehf} to ``lift'' the result to all the resolutions above the minimals, proving that the bounds of Conjecture~\ref{mainconj} hold for every module with that Hilbert function. 

Here is an example of this type of application.

\begin{example}\label{firstliftex}\em{
Let $S=k[a,b,c]$, and let $L \subset S$ be the lexicographic ideal such that $S/L$ has Hilbert function $H=(1,3,6,9,9,6,2)$. Let $I$ be the ideal $(a^3,b^4,c^4,b^2c^2)$. Then $L$ and $I$ have the same Hilbert function. We present the graded free resolutons of $S/L$ and $S/I$ using the Betti diagram notation of Macaulay 2. The columns and rows are numbered starting with zero, and $\beta_{ij}$ appears in row $j-i$ and column $i$ in the diagram. The Betti diagrams of $S/L$ and $S/I$ are below:

\bigskip

\begin{tabular}{ccccccccccccc}
$S/L$: & total: & 1 & 16 & 27 & 12 & \, & $S/I$: & total: & 1 & 4 & 5 & 2\\
& 0: & 1 & . & . & . & \, & & 0: & 1 & . & . & .\\
& 1: & . & . & . & . & \, & & 1: & . & . & . & .\\
& 2: & . & 1 & . & . & \, & & 2: & . & 1 & . & .\\
& 3: & . & 3 & 5 & 2 & \, & & 3: & . & 3 & . & .\\
& 4: & . & 5 & 9 & 4 & \, & & 4: & . & . & 2 & .\\
& 5: & . & 5 & 9 & 4 & \, & & 5: & . & . & 3 & .\\
& 6: & . & 2 & 4 & 2 & \, & & 6: & . & . & . & 2\\
\end{tabular}

\bigskip

Note that by making all potentially possible cancellations in the Betti diagram of $S/L$, we obtain the Betti diagram of $S/I$. Therefore $\beta^{S/I} \le$ $\beta^{S/J}$ for all ideals $J \subset S$ with the same Hilbert function as $I$; the resolution of $S/I$ is the unique minimal element in the partial order on resolutions with the fixed Hilbert function. The bounds on $S/I$ from Conjecture~\ref{mainconj} are \[27=\frac{1}{3!}(3)(6)(9) \le e(S/I)=36 \le \frac{1}{3!}(4)(7)(9)=42, \] and thus $S/I$ satisfies the conjecture. (This is actually immediate since $S/I$ has a quasipure resolution.) By Proposition~\ref{multsamehf}, the bounds of Conjecture~\ref{mainconj} hold for all modules with Hilbert function $H$.

Of course, not all examples will be this nice; there will often be incomparable minimal elements for a particular Hilbert function as a result of having the choice among cancellations, and then needs, at least \emph{a priori}, to check all minimals. (See Section~\ref{applications} for a way to avoid checking multiple configurations of graded Betti numbers in some cases.)
}
\end{example}

\begin{remark}
\label{notexist}\em{
While the Betti diagram following the cancellations in Example~\ref{firstliftex} is really the Betti diagram of a module, this need not be the case in applying the idea of Proposition~\ref{multsamehf}; the result is a purely numerical condition. Take all the potential Betti diagrams obtained from making as many consecutive cancellations in the resolution of a lexicographic ideal $L$ as numerically possible, and compute the conjectured bounds for each of these configurations. If all the potential Betti diagrams satisfy the bounds, we can conclude that all resolutions with the same Hilbert function as $L$ satisfy Conjecture~\ref{mainconj}.}
\end{remark}

We can use the same ideas to get a slight generalization of Theorem~\ref{pureres}, the result of Huneke and Miller on the pure resolution case, in codimension three. First, we give a remark we will use several times in the paper.

\begin{remark}
\label{cmartinian}\em{
Let $I$ be a homogeneous ideal of codimension $c$ in $R=k[x_1,\dots,x_n]$ such that $R/I$ is Cohen-Macaulay. Then there exists a homogeneous ideal $I' \subset R'=k[x_1,\dots,x_c]$ such that $R'/I'$ is Artinian and has the same graded Betti numbers (and thus multiplicity) as $R/I$. Hence we can reduce the Cohen-Macaulay case of Conjecture~\ref{mainconj} to the Artinian case.}
\end{remark}

\begin{proposition}
\label{codim3}
Let $R=k[x_1,\dots,x_n]$, and let $I$ be a homogeneous ideal in $R$ of codimension three such that $R/I$ is Cohen-Macaulay and has a pure resolution. Let $J$ be any homogeneous ideal in $R$ with the same Hilbert function such that $R/J$ is Cohen-Macaulay. Then Conjecture~\ref{mainconj} holds for $R/J$.
\end{proposition}

\begin{proof}
By Remark~\ref{cmartinian}, we may assume that $I$ and $J$ are Artinian. Theorem~\ref{pureres} implies that $R/I$ satisfies the bounds of Conjecture~\ref{mainconj}. We show that $\beta^{R/I} \le \beta^{R/J}$. 

Since the resolution of $R/I$ is pure, $I$ has minimal generators in only one degree, the lowest degree $d$ for which $H(R/I,d) < {n-1+d \choose d}$. Since $J$ has the same Hilbert function, $J$ has the same number of generators in degree $d$ and possibly more in higher degree. Hence $\beta^{R/I}_{1j}$ $\le \beta^{R/J}_{1j}$ for all $j$.

Because $I$ is Artinian, there is some degree $s$ in which $H(R/I)$ is last nonzero, and $\beta^{R/J}_{3,s+3}$ is the same for all ideals $J$ with the same Hilbert function as $I$. Also, since $R/I$ has a pure resolution, $\beta^{R/I}_{3t}=0$ for all $t \not = s+3$. Thus $\beta^{R/I}_{3j}$ $\le \beta^{R/J}_{3j}$ for all $j$.

Next we consider the second place in the resolution. Suppose $\beta^{R/I}_{2u} \not = 0$; then all other $\beta^{R/I}_{2v} = 0$. Let $L$ be the lexicographic ideal of $R$ with the same Hilbert function as $I$. Then \[ \beta^{R/I}_{2u} = \beta^{R/L}_{2u} - (\beta^{R/L}_{1u} + \beta^{R/L}_{3u}) \] since the resolution of $R/I$ is pure. Because the resolution of $R/J$ is also obtained from making cancellations in the resolution of $R/L$, we have \[ \beta^{R/I}_{2u} = \beta^{R/L}_{2u} - (\beta^{R/L}_{1u} + \beta^{R/L}_{3u}) \le \beta^{R/J}_{2u}. \] Therefore $\beta^{R/I} \le \beta^{R/J},$ and by Proposition~\ref{multsamehf}, Conjecture~\ref{mainconj} holds for $R/J$. \hfill $\square$
\end{proof}


\section{Computational work}\label{applications}

We discuss some additional ways to use Proposition~\ref{multsamehf} in this section. The results we obtain give some insight into why Conjecture~\ref{mainconj} is so difficult to prove in general, though to verify the conjecture for a particular Hilbert function, our techniques will often work reasonably well.

We used the computer algebra system Macaulay 2  to investigate Conjecture~\ref{mainconj} using Proposition~\ref{multsamehf}. Our goal was to use Proposition~\ref{multsamehf} for Artinian ideals to verify the upper bound in Conjecture~\ref{mainconj} for many Hilbert functions with small socle degree (the highest degree in which the Hilbert function is nonzero). We did most of our computations for ideals in $S=k[x_1,x_2,x_3]$.

To prove the upper bound of Conjecture~\ref{mainconj} for a particular Hilbert function, we began by computing the lexicographic ideal in $S$ corresponding to that Hilbert function. The idea was to keep the multiplicity fixed but to lower the $M_i$, proving that the multiplicity is bounded above by a smaller bound. Starting with the Betti diagram of the lexicographic ideal, we made all potentially possible cancellations in the diagram to find minimal sets of graded Betti diagrams in the partial order. At this stage, we computed the upper bound from the resulting potential Betti diagram(s). If the upper bound of Conjecture~\ref{mainconj} held for all the minimal diagrams, that proved the bound for all modules with that Hilbert function. 

There are some complications to consider, however. First, it is often difficult to determine whether a potential cancellation in a Betti diagram can occur. We decided to disregard this concern in our initial tests and return to it at the end; we will discuss this in more detail later in the section. Second, there may be many incomparable minimal elements in the partial order for a given Hilbert function. To minimize the number of Betti diagrams we had to check, we made our cancellations in a particular way that we outline below, finding a specific minimal element in the partial order, and it sufficed to check that diagram to determine if the upper bound held.

In three variables, the process is easy. For each (Artinian) Hilbert function, we first resolve the lexicographic ideal corresponding to that Hilbert function. $M_3$ is always fixed, for it is simply the maximal degree of a socle element plus three. We lower $M_1$ as much as possible by making all potentially possible cancellations of entries in the generators column of the Betti diagram with first syzygies. Finally, we lower $M_2$, if possible, by cancelling first syzygy entries in the Betti diagram with second syzygies. This yields a Betti diagram that may or may not represent the graded free resolution of an actual $S$-module. Regardless, we can still evaluate whether the bounds of Conjecture~\ref{mainconj} hold, computing the multiplicity directly from the Hilbert function and reading the minimum and maximum shifts from the diagram with the cancellations. The process lowers the product $M_1M_2M_3$ as much as possible; see Theorem~\ref{onediagram} and the surrounding discussion for the verification (in any number of variables).

To illustrate this method, we give an example.

\begin{example}\label{cancelex}\em{
Consider the Hilbert function $(1,3,6,7,3,1)$. The lexicographic ideal corresponding to this Hilbert function in $S=k[x_1,x_2,x_3]$ has the following resolution:

\bigskip
\begin{tabular}{cccccc} 
$S/L$: & total: & 1 & 12 & 19 & 8 \\
& 0: & 1 & . & . & . \\
& 1: & . & . & . & . \\
& 2: & . & 3 & 3 & 1 \\
& 3: & . & 6 & 10 & 4\\
& 4: & . & 2 & 4 & 2 \\
& 5: & . & 1 & 2 & 1\\
\end{tabular}

\bigskip

We first make all the cancellations potentially possible in the generators and first syzygies columns. This yields the Betti diagram below:

\bigskip
\begin{tabular}{cccccc} 
& total: & 1 & 6 & 13 & 8 \\
& 0: & 1 & . & . & . \\
& 1: & . & . & . & . \\
& 2: & . & 3 & . & 1 \\
& 3: & . & 3 & 8 & 4\\
& 4: & . & . & 3 & 2 \\
& 5: & . & . & 2 & 1\\
\end{tabular}

\bigskip

We then cancel any remaining pairs in the first and second syzygies columns, and obtain the following diagram:

\bigskip
\begin{tabular}{cccccc} 
& total: & 1 & 6 & 7 & 2 \\
& 0: & 1 & . & . & . \\
& 1: & . & . & . & . \\
& 2: & . & 3 & . & . \\
& 3: & . & 3 & 7 & 1\\
& 4: & . & . & . & . \\
& 5: & . & . & . & 1\\
\end{tabular}

\bigskip

The multiplicity corresponding to this Hilbert function is $1+3+6+7+3+1=21$. The maximum shifts are 4, 5, and 8. Because \[ 21 \le \frac{(4)(5)(8)}{3!} = \frac{80}{3}, \] the upper bound of Conjecture~\ref{mainconj} holds for the diagram with all the cancellations. But we have lowered the maximum shifts as much as possible, and therefore the upper bound holds for all modules $S/I$ with Hilbert function $(1,3,6,7,3,1)$.
}
\end{example}

Our computations were fruitful. In Macaulay 2, we tested all Hilbert functions for modules $S/I$ such that $I$ contains no linear forms and $S/I$ is zero in degree ten and higher; that is, we investigated the Hilbert functions of the form $(1,3,h_2,\dots,h_9)$. We tested over 677,000 Hilbert functions with this technique, and only 197 returned potential Betti diagrams that did not satisfy the upper bound of Conjecture~\ref{mainconj}. For each of those 197 examples, we looked at the potential Betti diagram with the cancellations more carefully, and we found reasons the diagrams could not represent the resolution of a module. Many cases were obvious; for example, an Artinian ideal in three variables that is not a complete intersection must have at least four generators. To study the remaining potential counterexamples, we used a theorem of Evans and Richert from \cite{EvansRichert}.

\begin{theorem}
\label{GrahamBen}
\em{(Evans-Richert)} Let $R=k[x_1, \dots, x_n]$, and let $M$ be a graded $R$-module. Let $i \ge 2$ be a positive integer. If $t$ is the smallest integer such that $\beta_{it}^M \not = 0$, then the sum of the graded Betti numbers $\beta_{i-1,j}^M$ for all $j<t$ is at least $i$.
\end{theorem}

One can think of this theorem as generalizing the fact that one needs at least two minimal generators of degree at least $d$ before having a minimal first syzygy of degree $d+1$. For some applications, see \cite{EvansRichert} or Chapter 2 of \cite{thesis}.

With Theorem~\ref{GrahamBen}, we proved that all but nine of the remaining potential counterexamples cannot exist. For the other examples, we needed a few easy computations. We give two examples below to illustrate the process of eliminating the 197 possible counterexamples.

\begin{example}\label{easyex}\em{
Consider the Hilbert function $(1,3,6,10,15,15,11)$ for quotients of $S=k[x_1,x_2,x_3]$. The resolution of the lexicographic ideal corresponding to this Hilbert function is below:

\bigskip

\begin{tabular}{cccccc} 
$S/L$: & total: & 1 & 25 & 43 & 19 \\
& 0: & 1 & . & . & . \\
& 1: & . & . & . & . \\
& 2: & . & . & . & . \\
& 3: & . & . & . & .\\
& 4: & . & 6 & 8 & 3 \\
& 5: & . & 7 & 12 & 5\\
& 6: & . & 12 & 23 & 11\\
\end{tabular}

\bigskip

After making all the cancellations in the appropriate order, we obtain the following Betti diagram:

\bigskip

\begin{tabular}{cccccc} 
$S/L$: & total: & 1 & 6 & 19 & 14 \\
& 0: & 1 & . & . & . \\
& 1: & . & . & . & . \\
& 2: & . & . & . & . \\
& 3: & . & . & . & .\\
& 4: & . & 6 & 1 & 3 \\
& 5: & . & . & . & .\\
& 6: & . & . & 18 & 11\\
\end{tabular}

\bigskip

The multiplicity corresponding to the Hilbert function is $1+3+6+10+15+15+11$ $=61$. The maximum shifts are 5, 8, and 9; note that the cancellations in degree eight are irrelevant to our computations. Since \[ 61 > \frac{1}{3!}(5)(8)(9)=60, \]the upper bound of Conjecture~\ref{mainconj} fails. We had to make the cancellations in degree six and seven to get maximum shifts for which the bound fails, so to show that the upper bound holds for all modules with this Hilbert function, it is enough to show that row four of this Betti diagram cannot exist. This is clear from Theorem~\ref{GrahamBen} (or by simply observing that there cannot be three minimal second syzygies of degree seven when there is a single first syzygy of degree at most six).
}
\end{example}

We demonstrate the technique used to eliminate the final nine potential counterexamples in the following example.

\begin{example}\label{hardex}\em{
Let $H$ be the Hilbert function $(1,3,6,10,15,17,17,17,15,10)$, which yields a multiplicity of 111. The lexicographic ideal in $S=k[x_1,x_2,x_3]$ corresponding to this Hilbert function has the following resolution:

\bigskip

\begin{tabular}{cccccc} 
$S/L$: & total: & 1 & 30 & 53 & 24 \\
& 0: & 1 & . & . & . \\
& 1: & . & . & . & . \\
& 2: & . & . & . & . \\
& 3: & . & . & . & .\\
& 4: & . & 4 & 4 & 1 \\
& 5: & . & 3 & 5 & 2\\
& 6: & . & 2 & 4 & 2\\
& 7: & . & 4 & 7 & 3 \\
& 8: & . & 6 & 12 & 6\\
& 9: & . & 11 & 21 & 10\\
\end{tabular}

\bigskip

We first make all the potential cancellations between the generators and first syzygies columns, and then we note that it will not be possible to lower $M_2$ (or $M_3$, of course) since we have 21 first syzygies of degree 11 but only six second syzygies in that degree. The diagram below is the result of these cancellations.

\bigskip

\begin{tabular}{cccccc} 
$S/L$: & total: & 1 & 4 & 27 & 24 \\
& 0: & 1 & . & . & . \\
& 1: & . & . & . & . \\
& 2: & . & . & . & . \\
& 3: & . & . & . & .\\
& 4: & . & 4 & 1 & 1 \\
& 5: & . & . & 3 & 2\\
& 6: & . & . & . & 2\\
& 7: & . & . & 1 & 3 \\
& 8: & . & . & 1 & 6\\
& 9: & . & . & 21 & 10\\
\end{tabular}

\bigskip

This gives maximum shifts of 5, 11, and 12, and the upper bound is thus 110, which is less than the multiplicity of 111. This Betti diagram cannot exist because of the single first syzygy of degree at most six and the second syzygy of degree seven. However, if we cancel the second syzygy of degree seven and a first syzygy of the same degree, Theorem~\ref{GrahamBen} does not forbid the resulting Betti diagram. Therefore we need to be a bit more creative.

We will show that an almost complete intersection $I$ of four degree five polynomials cannot have Hilbert function $H$. Call the generators of $I$ $f_1$, $f_2$, $f_3$, and $f_4$, and assume, reindexing if necessary, that the first three form a regular sequence. Let $J=(f_1,f_2,f_3):(f_4)$. We have the following short exact sequence: \[ 0 \longrightarrow S/J(-5) \longrightarrow S/(f_1,f_2,f_3) \longrightarrow S/I \longrightarrow 0. \]

We are given the Hilbert function of $S/I$, and the Hilbert function of the complete intersection $S/(f_1,f_2,f_3)$ is easy to compute. The Hilbert function of $S/J(-5)$ is the difference between the two:

\bigskip

\begin{tabular}{ccccccccccccccccc}
  & (1,&3,&6,&10,&15,&18,&19,&18,&15,&10,&6,&3,&1)\\
- & (1,&3,&6,&10,&15,&17,&17,&17,&15,&10)\\
\hline
  & (0,&0,&0,& 0,& 0,& 1,& 2,& 1,& 0,& 0,&6,&3,&1)\\
\end{tabular}

\bigskip

Therefore the Hilbert function of $S/J$ should be $(1,2,1,0,0,6,3,1)$, which is clearly not an allowable Hilbert function. Hence no almost complete intersection of four degree five polynomials exists with Hilbert function $H$, and thus the upper bound holds for all modules $S/I$ with this Hilbert function.
}
\end{example}

After ruling out the 197 possible counterexamples, we have the following theorem from the computations in Macaulay 2.

\begin{theorem}
\label{compub}
Let $I$ be a homogeneous ideal in $S=k[x_1,x_2,x_3]$ such that $S/I$ is zero in degree ten and higher. Then $S/I$ satisfies the upper bound of Conjecture~\ref{mainconj}; that is, \[ e(S/I) \le \frac{1}{3!}M_1M_2M_3.\]
\end{theorem}

\begin{remark}
\label{implications}\em{
This computational work in three variables has implications for possible proof techniques for Conjecture~\ref{mainconj}. First, the existence of the 197 potential counterexamples is significant. In the cases of Cohen-Macaulay modules $R/I$ with pure or quasipure resolutions, there are numerical proofs to show that the modules satisfy the bounds of Conjecture~\ref{mainconj}. Example~\ref{easyex} shows that such a numerical argument is impossible in general. Moreover, Example~\ref{easyex} proves that the quasipure result is the best possible in this direction; the Betti diagram in that example has only three nonzero rows above row zero, yet after cancellation, one gets a potential Betti diagram that does not satisfy the upper bound of Conjecture~\ref{mainconj}. In the other direction, the numerical counterexamples tell us that Conjecture~\ref{mainconj}, if true, gives extremely useful information about what resolutions can occur for a given Hilbert function. A positive solution to Conjecture~\ref{mainconj} would eliminate a number of candidates for minimal elements in the partial order on resolutions for a given Hilbert function that we currently cannot easily rule out.
}
\end{remark}

The process we used to prove Theorem~\ref{compub} is more complicated in more variables, but it is still possible to check only a single Betti diagram. To show this, we begin with a proposition that gives lower bounds for the maximal shifts in the resolution of a Cohen-Macaulay module. 

\begin{proposition}
\label{cmbetti}
Suppose $I$ is a homogeneous ideal in $R=k[x_1,\dots,x_n]$ of codimension $c$ such that $R/I$ is Cohen-Macaulay. Then in the resolution of $R/I$, $M_{i} \ge M_{i-1}+1$.
\end{proposition}

\begin{proof}
Suppose $R/I$ has regularity $d-c$. Since $R/I$ is Cohen-Macaulay, the dual of the resolution of $R/I$ is the minimal graded free resolution of the module $N=\text{Ext}^c_R(R/I,R)(-d)$. (The shift of $-d$ comes from the fact that the highest degree generator of a free module in the minimal free resolution of $R/I$ has degree $d$.) The Betti diagram of $M$ is obtained by rotating the Betti diagram of $R/I$ 180 degrees and shifting the degree.

From this, it is easy to see that $-m_{c-i}^N=M_i^{R/I}.$ Since $N$ is a graded module, $m_{c-i}^N+1 \le m_{c-i+1}^N$. Therefore \[-M_i^{R/I}+1 \le -M_{i-1}^{R/I}, \quad \hbox{and} \quad M_i^{R/I} \ge M_{i-1}^{R/I}+1.\] \hfill $\square$
\end{proof}

Here is a brief example to illustrate why the Cohen-Macaulay hypothesis is necessary in Proposition~\ref{cmbetti}.

\begin{example}\label{cmex}\em{
Consider the following resolution:

\bigskip

\begin{tabular}{cccccc} 
total: & 1 & 7 & 9 & 3 \\
0: & 1 & . & . & . \\
1: & . & . & . & . \\
2: & . & 6 & 8 & 3 \\
3: & . & 1 & 1 & .\\
\end{tabular}

\bigskip

This is the resolution of $k[a,b,c]/I$, where $I=(a^3,a^2b,a^2c,ab^2,abc,ac^2,b^4)$, so it is a quotient by a stable ideal, but the module is not Cohen-Macaulay. Note that there is a syzygy of degree five in column two but no syzygy of degree at least six in the third column. 

Juan Migliore kindly pointed out that $I$ is a basic double link, giving a nice way of constructing examples like this. (For a discussion of basic double linkage, see, e.g., \cite{migbook}.) Note that if $J$ is the ideal $(a^2,ab,ac,b^2,bc,c^2)$, then $I=aJ + (b^4)$. We could just as easily have let $I=aJ + (b^d)$, where $d \gg 0$, which allows $M_2-M_3$ to be as large as we want.
}
\end{example}

We now show that for each Artinian Hilbert function $H$, there exists an easily obtained potential Betti diagram $\mathcal D$ corresponding to $H$ such that if $\mathcal D$ satisfies the upper bound of Conjecture~\ref{mainconj}, then all modules with the same Hilbert function do as well. Our algorithm for producing $\mathcal D$ is the following: We lower $M_1$ as much as possible by making all potentially possible cancellations between the generator and first syzygy columns in the Betti diagram of the lexicographic ideal corresponding to $H$. Then we try to lower $M_2$ by cancelling first syzygies with second syzygies, and so on, proceeding left to right across the Betti diagram. Note that $M_n$ is fixed by the Hilbert function. (Also, there may be many such $\mathcal D$ corresponding to a single Hilbert function with the same product $M_1 \cdots M_n$.)

\begin{theorem}
\label{onediagram}
Let $H$ be a Hilbert function for an Artinian module of the form $R/I$. Let $L \subset R$ be the lexicographic ideal such that $H(R/L)=H$, and let $\mathcal D_L$ be the Betti diagram of $R/L$. Then the Betti diagram $\mathcal D$ obtained from consecutive cancellations in $\mathcal D_L$ with the above algorithm satisfies the property that if $J$ is any ideal with $H(R/J)=H$, then the product of the maximal shifts $M_1 \cdots M_n$ of $\mathcal D$ is at most $M_1^{R/J} \cdots M_n^{R/J}$.
\end{theorem}

\begin{proof}
The question is whether the algorithm gives us the lowest potentially possible value of $M_1M_2 \cdots M_n$ for a module with Hilbert function $H$. That is, we need to show that no other choice of consecutive cancellations in the resolution of the lexicographic ideal that would satisfy Proposition~\ref{cmbetti} gives a lower $M_1M_2 \cdots M_n$. Suppose we have a choice of cancellation at some point in the process. Consider the portion of a Betti diagram of a lexicographic ideal shown below, with columns $i-1$ through $i+1$ displayed. Suppose that the entries  $a,b,c$ are nonzero in degree $d$. By Proposition~\ref{cmbetti}, there are entries in row $d-i+1$ or below, in columns $i$ and $i+1$, that are nonzero. The asterisks represent possibly nonzero entries, and the column numbers are bolded at the top of the diagram.

\bigskip

\begin{tabular}{cccccccc}
& \boldmath{$0$} & \boldmath{$\dots$} & \boldmath{$i-1$} & \boldmath{$i$} & \boldmath{$i+1$} & \boldmath{$\dots$} \\
\hline
$d-i-1$: & . & \dots & * & * & $c$ & \dots \\
$d-i$:   & . & \dots & * & $b$ & * & \dots \\
$d-i+1$: & . & \dots & $a$ & $*$ & $*$ & \dots \\
\end{tabular}

\bigskip

We have a choice: We could cancel using the $a$ and $b$ in columns $i-1$ and $i$, attempting to lower $M_{i-1}$ and/or $M_i$, or we could cancel using the $b$ and $c$ in columns $i$ and $i+1$, attempting to lower $M_i$ and/or $M_{i+1}$. We claim that we must eventually make the cancellations in columns $i-1$ and $i$ to have any hope of lowering $M_i$ below $d+1$.

If we do not cancel the entry $a$ down to zero, we have $M_{i-1} \ge d$, and then Proposition~\ref{cmbetti} implies that $M_i \ge d+1$ and $M_{i+1} \ge d+2$. Therefore cancelling the entries $b$ and $c$ will not change $M_i$ and $M_{i+1}$.

Hence making all possible cancellations in the first two columns, then the second and third columns, etc., will lead to the a value of $M_1M_2 \cdots M_n$ less than or equal to the product of the maximal shifts of any ideal with Hilbert function $H$. \hfill $\square$
\end{proof}

\begin{remark}
\label{hardminimals}\em{
Our computational work on Conjecture~\ref{mainconj} helps to illustrate one reason that general progress has been so hard to obtain. When we used our algorithm in four variables, we found examples of potential Betti diagrams that do not satisfy the upper bound of Conjecture~\ref{mainconj}, and we were not able to show that all cannot occur. The Betti diagrams get much more complicated as one adds variables, and the dearth of results like Theorem~\ref{GrahamBen} probably precludes further significant progress with this technique right now. The main problem is that one cannot use this approach to get general results without being able to tell what resolutions actually occur at the bottom of the partial order on the resolutions of modules with a given Hilbert function. Conjecture~\ref{mainconj} has been so difficult at least in part because it is closely related to the problem of finding exactly what resolutions are the bottom of the posets, something we are a long way from being able to do. However, the technique of reducing to potential minimal elements in the partial order still gives a fast, easily tested sufficient condition for all modules with a particular Hilbert function to satisfy the conjecture(s).
}
\end{remark}

\section{Truncation}\label{truncation}

In this section, we discuss another way to investigate the upper bound of Conjecture~\ref{mainconj}. Throughout, $I$ will be a homogeneous ideal in $R=k[x_1, \dots, x_n]$ such that $R/I$ is Cohen-Macaulay. 

The technique we will use is to ``truncate'' $I$ in an appropriate degree to get rid of parts of the resolution in low degree that are irrelevant to the upper bound conjecture. Let $I_{\ge d}$ be the ideal in $R$ consisting of all elements of $I$ of degree $d$ or higher. Instead of working with $R/I$, we will work with modules of the form $R/I_{\ge d}$.

\begin{example}\label{truncate}\em{
Let $I=(a^3,b^3,c^3,ab,bc) \subset S=k[a,b,c]$. Truncating in degree three, we have $I_{\ge 3}=(a^3,b^3,c^3,a^2b,abc,ab^2,b^2c,bc^2)$. Note that $I$ contains two monomials of degree two that $I_{\ge 3}$ does not, and therefore, since both ideals are Artinian, we have $e(S/I)+2=e(S/I_{\ge 3})=13$. 
}
\end{example}

In Example~\ref{truncate}, the multiplicity increased when we truncated $I$. The next lemma shows how the multiplicity changes in general after truncation.

\begin{lemma}\label{truncatemult}
Let $I$ be a homogeneous ideal in $R=k[x_1,\dots,x_n]$. Then $e(R/I) \le e(R/I_{\ge d})$ for all $d$.
\end{lemma}

\begin{proof}
Fix $d$. Suppose first that $R/I$ is not Artinian. Then $R/I$ has a nonzero Hilbert polynomial. Since truncation only affects the dimension in finitely many degrees (that is, only finitely many fewer monomials are in $in_> (I_{\ge d})$ than are in $in_> (I)$), $R/I$ and $R/I_{\ge d}$ have the same Hilbert polynomial. Thus they have the same multiplicity.

If $R/I$ is Artinian, then $e(R/I)=\dim_k (R/I)$, and $e(R/I_{\ge d})=$ $\dim_k (R/I_{\ge d})$. Since $\dim_k (R/I)_t$ $\le \dim_k (R/I_{\ge d})_t$ for all $t$, $e(R/I)$ $\le e(R/I_{\ge d})$. \hfill $\square$
\end{proof}

Thus the multiplicity will always increase or stay the same after truncating. We are also interested in how the graded Betti numbers of $R/I$ are related to those of $R/I_{\ge d}$. We find this relation in the following lemma.

\begin{lemma}\label{truncatebetti}
Let $I$ be a homogeneous ideal in $R=k[x_1,\dots,x_n]$, and let $d$ be a positive integer. Then for each integer $l \ge 0$, \[ \beta_{i,i+d+l}^{R/I} = \beta_{i,i+d+l}^{R/I_{\ge d}}. \] That is, rows $d$ and higher of the Betti diagrams of $R/I$ and $R/I_{\ge d}$ are the same.

Moreover, if $I$ has its highest degree minimal generator in degree $\ge d$, and $R/I$ is Cohen-Macaulay, then $M_i^{R/I}=M_i^{R/I_{\ge d}}$ for all $1 \le i \le$ codim $I$.
\end{lemma}

\begin{proof}
We have the short exact sequence \[ 0 \longrightarrow I/I_{\ge d} \longrightarrow R/I_{\ge d} \longrightarrow R/I \longrightarrow 0. \] This induces a long exact sequence in Tor: For all $l \ge 0$, \[ \cdots \longrightarrow \text{Tor}_i(I/I_{\ge d},k)_{i+d+l} \longrightarrow \text{Tor}_i(R/I_{\ge d},k)_{i+d+l} \longrightarrow \text{Tor}_i(R/I,k)_{i+d+l} \longrightarrow \] \[ \text{Tor}_{i-1}(I/I_{\ge d},k)_{i+d+l} \longrightarrow \cdots \] is an exact sequence of $k$-vector spaces. Moreover, $I/I_{\ge d}$ has finite length; it is zero in degree $d$ and higher and has highest degree socle generator in degree $d-1$. Therefore $I/I_{\ge d}$ has regularity $d-1$, meaning $\beta^{I/I_{\ge d}}_{i,i+d+r}=0$ for all $i$ and all $r \ge 0$.

For $l \ge 0$, because $d+l>d-1$, \[ \dim_k \text{Tor}_i(I/I_{\ge d},k)_{i+d+l} = \beta_{i,i+d+l}^{I/I_{\ge d}} = 0. \] Similarly, \[ \dim_k \text{Tor}_{i-1}(I/I_{\ge d},k)_{i+d+l} = \beta_{i-1,i+d+l}^{I/I_{\ge d}} = 0. \] Consequently, as $k$-vector spaces, \[ \text{Tor}_i(R/I_{\ge d},k)_{i+d+l} \cong \text{Tor}_i(R/I,k)_{i+d+l}. \] Hence their dimensions over $k$ are equal, and thus for all $l \ge 0$, \[ \beta_{i,i+d+l}^{R/I} = \beta_{i,i+d+l}^{R/I_{\ge d}}. \] 

The final statement of the lemma follows immediately from the main portion of the lemma and Proposition~\ref{cmbetti}. \hfill $\square$
\end{proof}

\begin{example}\label{truncateres}\em{
We illustrate Lemma~\ref{truncatebetti} in an example. Let $S=k[a,b,c]$, and let \[I=(a^2,b^4,c^5,ac,bc,ab^2) \subset S.\] The Betti diagrams of $S/I$ and $S/I_{\ge 3}$ are below.

\bigskip

\begin{tabular}{ccccccccccccc}
$S/I$: & total: & 1 & 6 & 8 & 3 & \, & $S/I_{\ge 3}$: & total: & 1 & 10 & 15 & 6\\
& 0: & 1 & . & . & . & \, & & 0: & 1 & . & . & .\\
& 1: & . & 3 & 2 & . & \, & & 1: & . & . & . & .\\
& 2: & . & 1 & 2 & 1 & \, & & 2: & . & 8 & 11 & 4\\
& 3: & . & 1 & 2 & 1 & \, & & 3: & . & 1 & 2 & 1\\
& 4: & . & 1 & 2 & 1 & \, & & 4: & . & 1 & 2 & 1\\
\end{tabular}

\bigskip

We are truncating in degree three, and the Betti diagram in rows three and below does not change, just as we would expect from Lemma~\ref{truncatebetti}. Apart from the one in the (0,0) place, there are only zeros in the Betti diagram of $S/I_{\ge 3}$ before row two since we have no generators until degree three. The Betti diagram of $S/I_{\ge 3}$ has more generators (and syzygies) in row two than are in the Betti diagram of $S/I$ because all the monomials of degree three in $I$ are minimal generators of $I_{\ge 3}$.
}
\end{example}

We use Lemmas~\ref{truncatemult} and~\ref{truncatebetti} to reduce the upper bound portion of Conjecture~\ref{mainconj} to the case of ideals whose minimal generators are all in a single degree. We will only need the result in the Artinian case, so that is how we formulate it.

\begin{proposition}
\label{onedegree}
Let $I$ be an Artinian homogeneous ideal in $R=k[x_1,\dots,x_n]$. Let $d$ be the highest degree in which $I$ has a minimal generator. If $R/I_{\ge d}$ satisfies the upper bound of Conjecture~\ref{mainconj}, then so does $R/I$.
\end{proposition}

\begin{proof}
By Lemma~\ref{truncatemult}, $e(R/I) \le e(R/I_{\ge d})$, and we are assuming that $R/I_{\ge d}$ satisfies the upper bound of Conjecture~\ref{mainconj}, so \[ e(R/I) \le e(R/I_{\ge d}) \le \frac{1}{n!} \prod_{i=1}^n M^{R/I_{\ge d}}_i. \] By Lemma~\ref{truncatebetti}, $M^{R/I}_i = M^{R/I_{\ge d}}_i$ for each $i$, and thus $R/I$ satisfies the upper bound of Conjecture~\ref{mainconj}. \hfill $\square$
\end{proof}

As a consequence, we can reduce the upper bound question for Cohen-Macaulay ideals to a simpler case.

\begin{corollary}
\label{mcupreduction}
If the upper bound of Conjecture~\ref{mainconj} holds for all Artinian ideals in $R=k[x_1,\dots,x_n]$ whose minimal generators are all in a single degree, then the upper bound of Conjecture~\ref{mainconj} holds for all homogeneous ideals $I \subset R$ such that $R/I$ is Cohen-Macaulay.
\end{corollary}

\begin{proof}
We may immediately reduce to the Artinian case using Remark~\ref{cmartinian}. The reduction to ideals with minimal generators all in a single degree follows from Proposition~\ref{onedegree}: Given a homogeneous ideal $I$ with its highest degree generator in degree $d$, one can work with $R/I_{\ge d}$ instead of $R/I$. \hfill $\square$
\end{proof}

We use Corollary~\ref{mcupreduction} in conjunction with a result of Herzog and Srinivasan in \cite{HerzogSrin}. Recall that we say $R/I$ has a \textbf{quasipure resolution} if $M_{i-1}^{R/I} \le m_i^{R/I}$ for all $i > 1$. This means that the maximal shift at step $i-1$ is bounded above by the minimal shift at step $i$. 

Herzog and Srinivasan prove the following theorem in \cite{HerzogSrin}.

\begin{theorem}
\label{hsquasipure}
\emph{(Herzog-Srinivasan)} Let $I$ be a homogeneous ideal in $R$ such that $R/I$ is Cohen-Macaulay. If $R/I$ has a quasipure resolution, then $R/I$ satisfies both bounds of Conjecture~\ref{mainconj}.
\end{theorem}

Their proof is numerical: Any potential quasipure resolution below that of a lexicographic ideal satisfies the bounds; there is no need for there to exist a module with that resolution. We exploit this result by combining it with Proposition~\ref{cmbetti} and Lemmas~\ref{truncatemult} and~\ref{truncatebetti}.

\begin{theorem}
\label{highdegree}
Let $I$ be a homogeneous ideal of codimension $c$ in $R$ such that $R/I$ is Cohen-Macaulay of regularity $d$. Suppose $I$ has a minimal generator of degree $d$ or $d+1$. Then $R/I$ satisfies the upper bound of Conjecture~\ref{mainconj}; that is, \[ e(R/I) \le \frac{1}{c!} \prod_{i=1}^c M_i. \]
\end{theorem}

\begin{proof}
We may again assume that $I$ is Artinian and that $c=n$. Note that if $R/I$ has regularity $d$, then degree $d+1$ is the highest degree in which $I$ can have a minimal generator. By Proposition~\ref{cmbetti}, the resolution of $R/I_{\ge d}$ is concentrated in rows $d-1$ and $d$ of the Betti diagram, the bottom two rows. Therefore $R/I_{\ge d}$ has a quasipure resolution, and it satisfies the bounds of Conjecture~\ref{mainconj}. By Lemmas~\ref{truncatemult} and~\ref{truncatebetti}, we have \[ e(R/I) \le e(R/I_{\ge d}) \le \prod_{i=1}^c M_i^{R/I_{\ge d}} = \prod_{i=1}^c M_i^{R/I}. \] Hence $R/I$ satisfies the upper bound of Conjecture~\ref{mainconj}. \hfill $\square$
\end{proof}

This result gives, for example, an easy proof of the upper bound for stable Cohen-Macaulay ideals; Herzog and Srinivasan prove the upper bound for stable ideals without the Cohen-Macaulay hypothesis in \cite{HerzogSrin}. The examples in Section~\ref{applications} of potential Betti diagrams with three nonzero rows (instead of just two) that do not satisfy the upper bound of Conjecture~\ref{mainconj} show that no further reduction like the one in Theorem~\ref{highdegree} that relies on a result with a numerical proof is possible.

We illustrate Theorem~\ref{highdegree} with an example.

\begin{example}\label{highdeg}\em{
Let $I=(a^3,b^4,c^4,ab^2,a^2bc^3) \subset$ $S=k[a,b,c]$. Then $S/I$ has regularity six, and $I$ has a minimal generator in degree six. We resolve $S/I$ and $S/I_{\ge 6}$ below.

\bigskip

\begin{tabular}{ccccccccccccc}
$S/I$: & total: & 1 & 5 & 8 & 4 & \, & $S/I_{\ge 6}$: & total: & 1 & 27 & 46 & 20\\
& 0: & 1 & . & . & . & \, & & 0: & 1 & . & . & .\\
& 1: & . & . & . & . & \, & & 1: & . & . & . & .\\
& 2: & . & 2 & . & . & \, & & 2: & . & . & . & .\\
& 3: & . & 2 & 2 & . & \, & & 3: & . & . & . & .\\
& 4: & . & . & . & . & \, & & 4: & . & . & . & .\\
& 5: & . & 1 & 5 & 3 & \, & & 5: & . & 27 & 45 & 19\\
& 6: & . & . & 1 & 1 & \, & & 6: & . & . & 1 & 1\\
\end{tabular}

\bigskip

Note that $S/I$ does not have a quasipure resolution, so Theorem~\ref{hsquasipure} does not apply. However, the truncation $S/I_{\ge 6}$ does have a quasipure resolution, and $S/I$ and $S/I_{\ge 6}$ have the same maximum shifts at each step in the resolution. Also, $e(S/I) = 31 \le$ $57 = e(S/I_{\ge 6})$, and thus $S/I$ satisfies the upper bound of Conjecture~\ref{mainconj} because $S/I_{\ge 6}$ does.
}
\end{example}

It would be interesting to have structure theorems for Artinian ideals whose minimal generators are all in one degree, the situation in Corollary~\ref{mcupreduction}. For example, can we say anything about their Hilbert functions, perhaps giving more detailed upper bounds on their growth than what follows from Macaulay's Theorem? Even more ambitiously, what can we say about the graded free resolutions of such ideals? These questions are difficult, but perhaps it would be possible to make some progress in the codimension three case. The upper bound of Conjecture~\ref{mainconj} is wide open even in the Cohen-Macaulay codimension three case, so such results would represent substantial progress.

\end{document}